# О геометрических медианах треугольников[*]

Петр Панов

Национальный исследовательский университет
«Высшая школа экономики», panovpeter@mail.ru

25 января 2021 г.

**Аннотация**. Получены универсальные оценки аффинного типа для возможного расположения геометрических медиан периметров треугольников и для расположения геометрических медиан треугольных областей. В конце обсуждаются некоторые альтернативные реализации пространства треугольников.

**Ключевые слова**. Геометрическая медиана, пространство треугольников, вырожденный треугольник, медианное отображение

**Введение**

Геометрическая медиана подмножества $S$ метрического пространства $M$ – это точка $m \in M$, суммарное расстояние от которой до всех точек множества $S$ минимально. В простейшем случае, когда $S$ является конечным подмножеством числовой прямой $R$ или числовой выборкой, это определение можно записать в виде

$$m = \underset{X \in R}{\arg\min} \sum_{X_i \in S} |X_i - X|,$$

и оно фактически совпадает с определением статистической медианы выборки $S$. Так что геометрическая медиана является естественным пространственным обобщением статистической медианы.

С каждым треугольником, расположенным в евклидовой плоскости $\mathbb{R}^2$, можно связать три сорта геометрических медиан. Во-первых, это геометрическая медиана $m_0$ трех вершин треугольника $A, B$ и $C$

$$m_0 = \underset{X \in R^2}{\arg\min}(|A - X| + |B - X| + |C - X|). \tag{1}$$

Далее, это геометрическая медиана $m_1$ его периметра, состоящего из отрезков $a, b, c$,

$$m_1 = \underset{X \in R^2}{\arg\min} \left( \int_{P \in a} |P - X| \, dP + \int_{P \in b} |P - X| \, dP + \int_{P \in c} |P - X| \, dP \right). \tag{2}$$

---

[*] Выражаю признательность А.В. Савватееву за внимание, проявленное к работе

Здесь суммарное расстояние от точки $X$ до сторон треугольника вычисляется уже с помощью интегралов. Наконец, $m_2$ – это геометрическая медиана внутренности треугольника – треугольной области, обозначенной здесь $abc$,

$$m_2 = \arg\min_{X \in \mathrm{R}^2} \int_{P \in abc} |P - X|\, dP. \qquad (3)$$

В этом случае суммарное расстояние вычисляется с помощью двойного интеграла.

На рисунке 1 показано расположение всех трех геометрических медиан треугольника $ABC$ со сторонами 9,7,5.

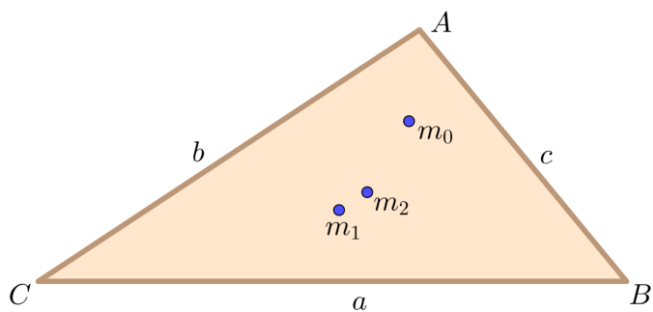

Рис. 1. Три типа геометрических медиан в треугольнике

Геометрическая медиана $m_0$ известна как точка Ферма – Торричелли треугольника, и существуют простые формулы для ее вычисления (Kimberling, 2020, точка X(13)).

Геометрические медианы $m_1, m_2$ в списке Кимберлинга отсутствуют, и это означает, что такие формулы для них не известны. Наша основная цель заключается в получении как можно более точной оценки для местоположения точек $m_1$ и $m_2$ внутри треугольника. Важный промежуточный шаг в этом направлении – это нахождение точных формул для геометрических медиан $m_1$ и $m_2$ вырожденных треугольников (Предложение 3).

**Основной результат**

Сформулируем основной результат настоящей работы в геометрической форме.

**Утверждение 1**. *Пусть задан произвольный треугольник $\delta$, и пусть*

- *$\delta'$ – это его образ при гомотетии с коэффициентом $1/4$ и с центром, расположенным в центроиде $\delta$*
- *$\delta''$ – это криволинейный треугольник, составленный из дуг трех гипербол, при этом вершины $\delta''$ совпадают с вершинами $\delta'$, а каждая из трех гипербол имеет центром одну из вершин $\delta$, касается двух медиан $\delta$ и имеет асимптотами две стороны $\delta$ (рис. 2)*

*Тогда*

1) *геометрическая медиана $m_1$ треугольника $\delta$ лежит внутри треугольника $\delta'$*

2) *геометрическая медиана $m_2$ треугольника $\delta$ лежит внутри треугольника $\delta''$.*

Все треугольники, присутствующие в Утверждении 1, изображены на рисунке 2, где дополнительно еще проведены медианы исходного треугольника $\delta$.

**Замечание 1.** Обе предъявленные здесь оценки – аффинного типа. Дело в том, что при аффинном преобразовании треугольника $\delta$ точно такому же преобразованию подвергаются связанные с ним треугольники $\delta'$ и $\delta''$.

Как будет показано дальше, обе полученные оценки в некотором смысле являются наилучшими. В несколько иной, аналитической форме, они будут представлены в Предложении 3.

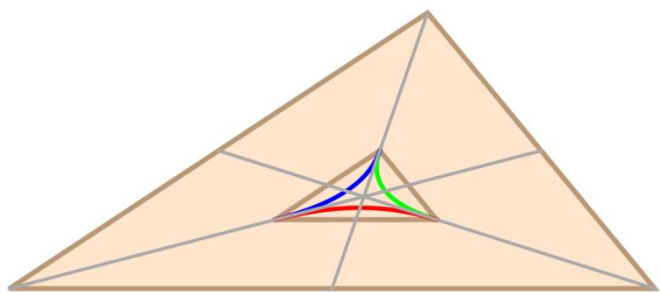

Рис. 2. Треугольники $\delta$, $\delta'$ и криволинейный треугольник $\delta''$, составленный из дуг гипербол

В дальнейшем изложении мы будем использовать результаты компьютерных экспериментов и некоторые аналитические вычисления, заимствованные из предыдущих работ (Панов, 2018; Панов, Савватеев, 2020).

**Пространство треугольников и медианные отображения $\mathcal{M}_0, \mathcal{M}_1, \mathcal{M}_2$**

Из определений (1) – (3) легко следует, что при преобразовании подобия каждая из медиан $m_i$ исходного треугольника переходит в соответствующую медиану преобразованного треугольника. А так как каждый треугольник подобен треугольнику с периметром 1, при исследовании расположения геометрических медиан только такими треугольниками можно и ограничиться. Пространство треугольников единичного периметра допускает множество геометрических реализаций (Behrend, 2014), используем одну из простейших.

Введем в рассмотрение правильный треугольник $\Delta$ с вершинами $\left(-\frac{1}{\sqrt{3}}, -\frac{1}{3}\right), \left(\frac{1}{\sqrt{3}}, -\frac{1}{3}\right)$ и $\left(0, \frac{2}{3}\right)$ – это большой треугольник на рисунке 3, его высота равна 1, а его центр расположен в начале координат.

Выберем произвольную точку внутри треугольника и опустим из нее перпендикуляры на левую, правую и нижнюю стороны треугольника. Обозначим длины перпендикуляров в этом порядке символами $a, b$ и $c$. Ясно, что сумма этих трех чисел равна высоте треугольника $\Delta$, то есть $a + b + c = 1$. Чтобы объявить $a, b$ и $c$ длинами сторон треугольника, недостает только неравен-ства треугольника. Нетрудно сообразить, что эти неравенства выполняются для точек внутреннего серединного треугольника $\nabla = -\frac{1}{2}\Delta$, заполненного цветными точками на рисунке 3. Таким образом, треугольник $\nabla$ на самом деле представляет пространство всех треугольников с периметром 1.

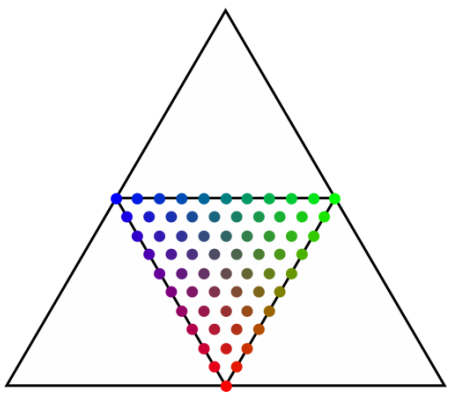

Рис 3. Высота большого треугольника равна 1; внутренний треугольник $\nabla$, заполненный цветными точками, – это пространство треугольников периметра 1

Граница этого пространства, которую мы будем обозначать $\partial\nabla$, соответствует вырожденным треугольникам, у которых одна из сторон равна сумме двух других (рис. 4).

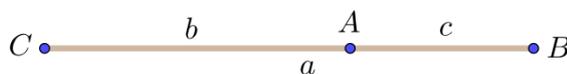

Рис. 4. Вырожденный треугольник, $a = b + c$

При построении медианных отображений мы будем использовать барицентрические координаты, связанные с треугольником. Напомним их определение. Треугольная область $abc$ – это выпуклая оболочка вершин треугольника $A, B$ и $C$. Значит, каждая точка $P \in abc$, представляется в виде

$$P = \lambda_a A + \lambda_b B + \lambda_c C,$$

где

$$\lambda_a + \lambda_b + \lambda_c = 1.$$

Эти три числа $\lambda_a, \lambda_b, \lambda_c$ и называются барицентрическими координатами точки $P$ относительно треугольника $ABC$ (Балк, Болтянский, 1987). Для внутренних точек треугольника барицентрические координаты положительны. На каждой из сторон треугольника одна из координат равна 0, а вершины треугольника имеют координаты $(1,0,0), (0,1,0), (0,0,1)$.

Перейдем к построению медианных отображений. Вместе с треугольником $\Delta$ рассмотрим три таких же правильных треугольника $\Delta_0, \Delta_1, \Delta_2$ с высотой, равной 1. Для каждого $i = 0,1,2$ опре-делим медианное отображение $\mathcal{M}_i: \nabla \to \Delta_i$ из пространства треугольников единичного периметра $\nabla$ в треугольник $\Delta_i$. Возьмем точку $x$, принадлежащую пространству $\nabla$. Ей соответствует некоторый треугольник $abc$. Вычислим барицентрические координаты его геометрической медианы $m_i$ – числа $\lambda_a^i, \lambda_b^i, \lambda_c^i$. По определению медианное отображение $\mathcal{M}_i$ ставит в соответствие точке $x \in \nabla$ точку с теми же самыми барицентрическими координатами $\lambda_a^i, \lambda_b^i, \lambda_c^i$ в треугольнике $\Delta_i$.

**Численные эксперименты: образы медианных отображений**

На рисунках 5–7 представлены образы медианных отображений $\mathcal{M}_0, \mathcal{M}_1, \mathcal{M}_2$, точнее образ достаточно плотной решетки пространства треугольников $\nabla$, типа той, что изображена на рисунке 3. На каждом из рисунков 5–7 точка образа наследует цвет точки прообраза из пространства $\nabla$.

Образ отображения $\mathcal{M}_0$, заполняет весь треугольник $\Delta_0$ (рис. 5). На координатном языке это означает, что для любых положительных чисел $\lambda_a, \lambda_b, \lambda_c$, сумма которых равна 1, существует треугольник, геометрическая медиана $m_0$ которого имеет эти числа своими барицентрическими координатами.

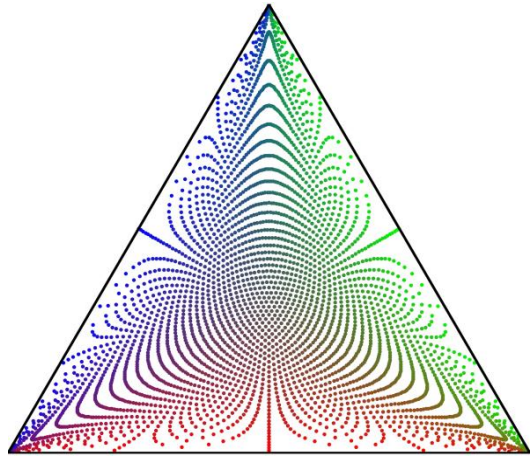

Рис. 5. Образ отображения $\mathcal{M}_0$ в треугольнике $\Delta_0$

Сравнение рисунков 3 и 5 показывает, что отображение $\mathcal{M}_0$ осуществляет частичное раздутие вершин треугольника $\nabla$ и стягивание его приграничных областей. Все это также легко следует из явных формул для вычисления барицентрических координат геометрической медианы $m_0$, то есть точки Ферма – Торричелли (Kimberling, 2020).

Перейдем к отображению $\mathcal{M}_1$. Из-за отсутствия явных формул для геометрической медианы $m_1$ периметра треугольника, для построения образа отображения $\mathcal{M}_1$, минимизационную задачу (2) приходится решать в полном объеме.

На рисунке 6 видно, что образ отображения заполняет треугольник $\Delta_1/4$, локализованный вблизи центра треугольника $\Delta_1$. Отображение не является однозначным, и чтобы выявить его структуру вблизи границы, нам пришлось деформировать решетку отображаемых точек и *дополнительно убрать из нее некоторый объем точек*.

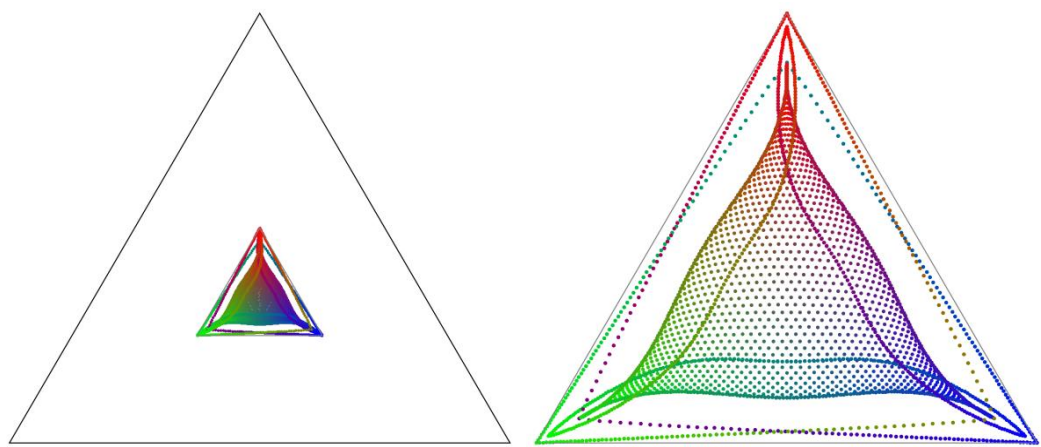

Рис. 6. Образ отображения $\mathcal{M}_1$ внутри треугольника $\Delta_1$ и его увеличенное изображение

Обращает на себя внимание, что при отображении $\mathcal{M}_1$ граница $\partial\nabla$ двукратно накрывает границу треугольника $\Delta_1/4$. В следующем разделе мы вернемся к обсуждению этого факта. А пока рекомендуем обратиться к анимации (Панов, 2020), где в динамике показано как формируется образ отображения $\mathcal{M}_1$, изображенный на рисунке 6.

Образ отображения $\mathcal{M}_2$ тоже содержится в треугольнике $\Delta_1/4$, но не заполняет его целиком. Он заполняет каспоидальный треугольник, который, как будет объяснено дальше, ограничен тремя отрезками гиперболы (рис. 7 справа).

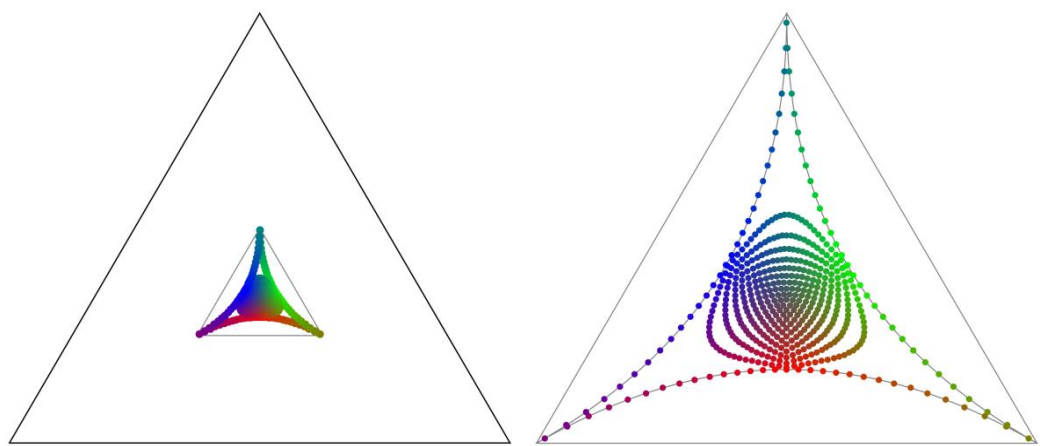

Рис. 7. Образ отображения $\mathcal{M}_2$ и его увеличенное изображение

Подведем первые итоги наших экспериментов. Во-первых, намного легче работать с теми отображениями, для которых известны простые аналитические формулы (имеется ввиду отображение $\mathcal{M}_0$). Во-вторых, как это видно из рисунков 5–7, граница образа каждого из отображений $\mathcal{M}_i$ формируется приграничными точками пространства $\nabla$. А как уже говорилось, эта граница $\partial\nabla$ со-стоит из вырожденных треугольников. Так что необходимо научиться вычислять барицентрические координаты геометрических медиан вырожденных треугольников. Этим мы сейчас и займемся.

**Продолжение отображений $\mathcal{M}_1$ и $\mathcal{M}_2$ на границу $\partial\nabla$**

Сложность заключается в том, что хотя определение (2) работает для вырожденного треугольника, и мы можем найти его геометрическую медиану $m_1$, но не совсем понятно, что это такое – барицентрические координаты точки относительно вырожденного треугольника. С геометрической медианой $m_2$ вырожденного треугольника дела обстоят хуже, для нее теряет смысл само определение (3). В обоих случаях при работе с геометрическими медианами вырожденных треугольников приходится использовать соображения непрерывности.

Полезно будет представлять себе вырожденный треугольник $abc$ как предел обычных треугольников $a_n b_n c_n$ таких, что $a_n \to a, b_n \to b, c_n \to c$. Коротко мы будем писать $a_n b_n c_n \to abc$. Сформулируем следующее предложение-определение, которое позволяет корректно определить барицентрических координат медиан $m_1$ и $m_2$ вырожденного треугольника.

**Предложение 1**. *Пусть $abc$ – это вырожденный треугольник. Рассмотрим последовательность сходящих к нему обычных треугольников $a_n b_n c_n \to abc$. Барицентрические координаты геометрической медианы $m_i, i = 1,2$ треугольника $a_n b_n c_n$ обозначим – $\lambda^i_{a_n}, \lambda^i_{b_n}, \lambda^i_{c_n}$. Тогда существуют пределы*

$$\lambda^i_{a_n} \to \lambda^i_a, \lambda^i_{b_n} \to \lambda^i_b, \lambda^i_{c_n} \to \lambda^i_c,$$

*которые не зависят от выбора последовательности $a_n b_n c_n$. Определенные таким образом числа $\lambda^1_a, \lambda^1_b, \lambda^1_c$ мы и будем называть барицентрическими координатами геометрической медианы $m_1$ треугольника $abc$, а числа $\lambda^2_a, \lambda^2_b, \lambda^2_c$ координатами его геометрической медианы $m_2$.*

**Замечание 2.** Предложение 1 говорит о возможности непрерывного продолжения медианных отображений $\mathcal{M}_1$ и $\mathcal{M}_2$ пространства треугольников $\nabla$ на его границу $\partial\nabla$, состоящую из вырожденных треугольников.

**Замечание 3.** Обратим внимание на то, что отображение $\mathcal{M}_0$ не допускает непрерывного продолжения на границу $\partial\nabla$. Об этом свидетельствует факт раздутия угловых точек пространства $\nabla$, который мы зафиксировали при обсуждении рисунка 5. Мы еще вернемся к обсуждению этого вопроса в последнем разделе статьи, где будут упомянуты другие реализации пространства треугольников, с другими граничными свойствами.

Дадим теперь конструктивное геометрическое истолкование ситуации, описанной в Предложении 1.

**Предложение 2**. *Пусть $abc$ – это вырожденный треугольник и $a \geq b \geq c$, и пусть $a_n b_n c_n$ – это последовательность сходящих к нему обычных треугольников. Тогда*

1) *геометрическая медиана $m_1$ треугольника $a_n b_n c_n$ удалена от общей вершины сторон $a_n$ и $b_n$ на расстояние, близкое к $a/2$, а отношение расстояний от нее до сторон $a_n$ и $b_n$ близко к 1*

2) *геометрическая медиана $m_2$ треугольника $a_n b_n c_n$ удалена от общей вершины сторон $a_n$ и $b_n$ на расстояние, близкое к $\sqrt{ab/2}$, и отношение расстояний от нее до сторон $a_n$ и $b_n$ близко к 1*

**Замечание 4.** На менее формальном и более выразительном языке, можно было бы сказать, что *геометрические медианы $m_1$ и $m_2$ вырожденного треугольника $abc$, $a \geq b \geq c$ лежат на биссектрисе его меньшего угла на расстояниях $a/2$ и $\sqrt{ab/2}$ от вершины этого угла*.

Последнее утверждение проиллюстрировано с помощью рисунка 8, на котором изображен треугольник, близкий к вырожденному треугольнику со сторонами $9, 8, 1$ и растянутый по вертикали в 1000 раз. При этом растяжении его биссектриса перешла в голубой наклонный отрезок. Рассчитанное для этого треугольника положение точек $m_1$ и $m_2$, действительно, близко́ к указанному в Предложении 2 и Замечании 3.

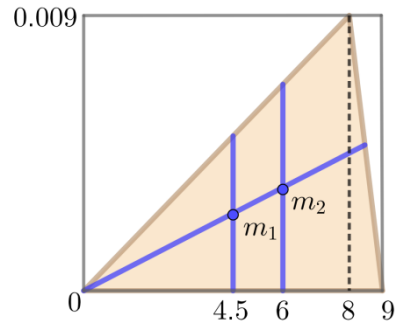

Рис. 8. Треугольник растянут вдоль вертикали в 1000 раз

Дадим теперь набросок доказательства Предложения 2.

Итак, $a$ – это бо́льшая сторона вырожденного треугольника $abc$ (рис. 4). В этом случае определение (2) геометрической медианы $m_1$ треугольника $abc$ можно переписать в виде

$$m_1 = \arg\min_{X \in \mathbb{R}^2} 2 \int_{P \in a} |P - X|\, dP.$$

Откуда сразу следует, что точку минимума нужно искать среди точек $X \in a$ и что на самом деле этот минимум достигается в середине отрезка $a$. Значит, геометрическая медиана $m_1$ – это середина бо́льшей стороны $a$. Тут нам не понадобился никакой предельный переход. Но в то же время отсюда легко вывести, что для обычных треугольников, близких к $abc$, геометрическая медиана $m_1$ расположена вблизи середины их бо́льшей стороны.

Этот же результат можно получить, приближенно решая градиентную систему для минимизационной задачи (2). Для обычных треугольников, близких к вырожденному треугольнику $abc$, в первом приближении опять получаем, что точка $m_1$ делит большую сторону треугольника пополам. А второе приближение как раз говорит о том, что расстояния от $m_1$ до двух больших стороны треугольника равны между собой.

Градиентная система для геометрической медианы $m_2$, то есть для задачи (3), выписана в работе (Панов, Савватеев, 2020). Одна из форм записи этой системы была использована там для непосредственного доказательства удаленности точки $m_2$ от вершины вырожденного треугольника $abc$, $a \geq b \geq c$ на расстояние $\sqrt{ab/2}$. Равенство расстояний от $m_2$ до двух больших стороны треугольника также следует из анализа приближенного решения градиентной системы для задачи (3).

**Точные формулы для геометрических медиан вырожденных треугольников**

Следующее утверждение является прямым переводом Предложения 2 на язык барицентрических координат.

**Предложение 3.** *Пусть $abc$ – это вырожденный треугольник, в котором $a \geq b \geq c$. Тогда*

1) *геометрическая медиана его периметра, точка $m_1$, имеет барицентрические координаты*

$$\lambda_a = \frac{a}{4b}, \qquad \lambda_b = \frac{1}{4}, \qquad \lambda_c = \frac{3}{4} - \frac{a}{4b}. \tag{4}$$

2) *геометрическая медиана его внутренности, точка $m_2$, имеет барицентрические координаты*

$$\lambda_a = \frac{1}{2\sqrt{2}}\sqrt{a/b}, \qquad \lambda_b = \frac{1}{2\sqrt{2}}\sqrt{b/a}, \qquad \lambda_c = 1 - \frac{1}{2\sqrt{2}}\left(\sqrt{a/b} + \sqrt{b/a}\right). \tag{5}$$

Отметим, что формулы (4) позволяют непрерывно продолжить отображение $\mathcal{M}_1$ на границу пространства $\nabla$ и подтвердить, что при этом $\partial \nabla$ дважды накрывает границу треугольника $\Delta_1/4$, как это изображено на рисунке 2.

В свою очередь формулы (5) позволяют непрерывно продолжить отображение $\mathcal{M}_2$ на границу пространства $\nabla$. Эта граница под действием отображения $\mathcal{M}_2$ переходит в криволинейный треугольник (рис. 7), сторонами которого являются отрезки гипербол, одна из которых имеет уравнение

$$x(t) = \frac{1}{\sqrt{6}}\operatorname{sh} t, \qquad y(t) = \frac{2}{3} - \frac{1}{\sqrt{2}}\operatorname{ch} t, \qquad -\ln\sqrt{2} < t < \ln\sqrt{2}.$$

После этого остается только сказать, что в силу аффинной эквивариантности барицентрических координат Утверждение 1 является прямым следствием Предложения 3.

**Градиентная система и интегралы**

В работе (Панов, Савватеев, 2020) показано, что градиентную систему для нахождения геометрической медианы $m_2$ ограниченной области $\Omega \subset \mathbb{R}^n$ можно записать в виде

$$\int_{P \in \partial \Omega} |P - X| \vec{n}(P) dP = 0,$$

где $\vec{n}(P)$ – единичный вектор внешней нормали к границе области $\partial \Omega$ в точке $P$. Для треугольной области $abc$ эта система трансформируется к виду

$$\left(\int_{P \in a} |P - m| dP\right) \vec{n}_a + \left(\int_{P \in b} |P - m| dP\right) \vec{n}_b + \left(\int_{P \in c} |P - m| dP\right) \vec{n}_c = 0$$

и далее к виду

$$\frac{1}{a} \int_{P \in a} |P - X| dP = \frac{1}{b} \int_{P \in b} |P - X| dP = \frac{1}{c} \int_{P \in c} |P - X| dP = 0.$$

Так или иначе, здесь, а также при работе с геометрической медианой $m_1$ (формула (2)), возникают интегралы типа $\int_{P \in a} |P - m| dP$. Если считать, что отрезок $a$ лежит на оси абсцисс – $a = [\alpha, \beta]$, точка $P$ имеет координаты $P = (\xi, 0)$, а точка $X$ имеет координаты $X = (x, y)$, то этот интеграл принимает вид $I(x, y) = \int_\alpha^\beta \sqrt{(\xi - x)^2 + y^2} \, d\xi$. И вот две формулы для его вычисления

$$I(x, y) = \frac{(\beta - x)\sqrt{(\xi - x)^2 + y^2} + (x - \alpha)\sqrt{(x - \alpha)^2 + y^2}}{2} + $$
$$+ \ln\left(\left(\sqrt{(\xi - x)^2 + y^2} + (\beta - x)\right)\left(\sqrt{(x - \alpha)^2 + y^2} + (x - \alpha)\right)\right) - \frac{y^2}{2} \ln y^2 \quad (6)$$

или

$$I(x, y) = \frac{(\beta - x)\sqrt{(\xi - x)^2 + y^2} + (x - \alpha)\sqrt{(x - \alpha)^2 + y^2}}{2} +$$
$$+ \frac{y^2}{2} \ln \frac{\sqrt{(x - \alpha)^2 + y^2} + (x - \alpha)}{\sqrt{(\beta - x)^2 + y^2} - (\beta - x)}. \quad (7)$$

Для треугольников, близких к вырожденным, можно считать, что величина $y$ в таких интегралах мала, так как обе геометрические медианы $m_1$ и $m_2$ лежат внутри треугольника (Панов, 2018). И при вычислении асимптотики интегралов в случае $x \in [\alpha, \beta]$ удобнее пользоваться формулой (6), а например когда $x > \beta$, удобнее использовать формулу (7).

Геометрические медианы $m_1$ и $m_2$, треугольника, близкого к вырожденному, лежат между двумя его бо́льшими сторонами. Формула (6) показывает, что для этих двух сторон при вычислении соответствующих им интегралов вида $I(x, y)$ появляются логарифмические слагаемые типа $y^2 \ln y^2$. В то же время формула (7) показывает, что для меньшей стороны логарифмическое слагаемое отсутствует. Именно это обстоятельство объясняет, почему в треугольниках, близких к вырожденным, расстояния от геометрической медианы $m_1$ до двух бо́льших сторон треугольника асимптотически равны между собой, и точно так же для $m_2$ (Предложение 2).

**Добавление: Медианное отображение $\mathcal{M}_0$, еще один тип вырожденных треугольников и другие реализации пространства треугольников**

Вернемся к медианному отображению $\mathcal{M}_0 \colon \nabla \to \Delta_0$, которое сопоставляет геометрической медиане $m_0$, то есть точке Ферма – Торричелли, набор ее барицентрических координат. Напомним, что нам не удалось продолжить это отображение на границу $\partial \nabla$. Причина заключалась в раздутии угловых точек пространства $\nabla$ при отображении $\mathcal{M}_0$ (рисунок 5 и Замечание 4). Эти угловые точки соответствуют вырожденным равнобедренным треугольникам с боковыми сторонами длины ½, и основанием нулевой длины. На следующем рисунке 9 представлено однопараметрическое семейство треугольников, близких к одному из таких вырожденных треугольников.

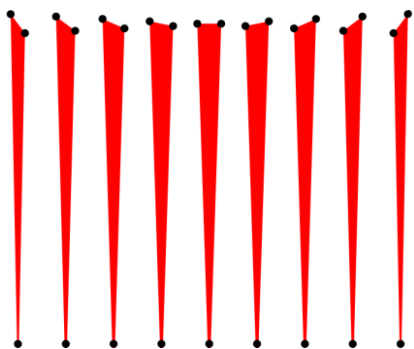

Рис. 9. Семейство треугольников, близких к равнобедренному вырожденному треугольнику со сторонами ½, ½, 0

Простые геометрические соображения, а также явные формулы для вычисления барицентрических координат точки Ферма – Торричелли (Kimberling, 2020, точка X(13)), показывают, что при медианном отображении $\mathcal{M}_0$ точки, соответствующие треугольникам из этого семейства, действительно, располагаются вдоль одной из сторон треугольника $\Delta_0$, что соответствует рисунку 5.

До сих пор мы работали с вырожденными треугольниками, у которых бо́льшая сторона равна сумме двух других (рис. 4), дальше мы их будем называть вырожденными треугольниками *первого типа*. Но предыдущие рассуждения показывают, что в рассмотрение следует еще ввести равнобедренные треугольники с нулевым углом при вершине и двумя другими углами при основании, сумма которых равна 180°. Такие треугольники мы будем называть вырожденными треугольниками *второго типа*. Это вырожденные треугольники, близки к тем, что изображены на рисунке 9.

При построении пространства треугольников $\nabla$ мы воспользовались тем, что любой треугольник однозначно определяется длинами своих сторон $a, b, c$. Но точно так же любой треугольник с точностью до подобия определяется величинами своих углов $\alpha, \beta, \gamma$. Сопоставим теперь такому треугольнику точку с барицентрическими координатами $\alpha/\pi$, $\beta/\pi$, $\gamma/\pi$ внутри фиксированного правильного треугольника $\Delta$. При таком сопоставлении треугольник $\Delta$ становится еще одной реализацией пространства треугольников. При этом граница $\partial\Delta$ оказывается составленной из вырожденных треугольников второго типа, и здесь уже медианное отображение $\mathcal{M}_0 : \Delta \to \Delta_0$ допускает непрерывное продолжение на границу $\partial\Delta$.

Как уже упоминалось, в работе (Behrend, 2014) рассмотрено большое количество геометрических реализаций пространства треугольников. Рассмотренное нами пространство $\nabla$ обозначено там символом $\mathcal{N}$, и там же представлено наглядное изображение этого пространства (Behrend, 2014, p. 34). Мы воспроизводим здесь эту иллюстрацию для пространства $\nabla$ с небольшими изменениями, а также аналогичную картинку для пространства $\Delta$. На рисунке 10 к точкам, равномерно заполняющим эти пространства, прикреплены треугольники, соответствующие этим точкам.

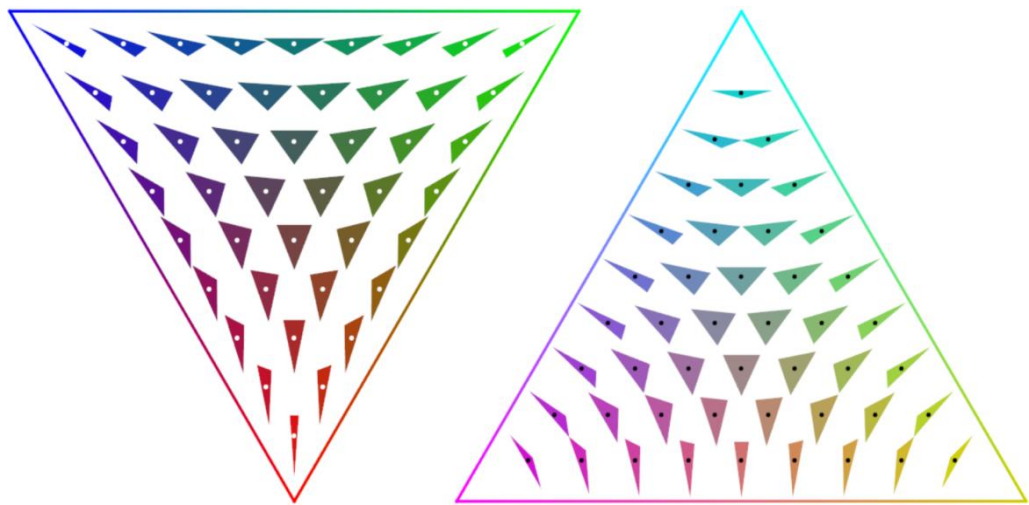

Рис. 10. Две реализации пространства треугольников – $\nabla$ и $\Delta$

На этом рисунке видно, что вырожденные треугольники первого типа сосредоточены вдоль границы пространства $\nabla$, а вырожденные треугольники второго типа – вдоль границы пространства $\Delta$. Имеется привлекательная возможность объединить эти два пространства в одно целое, что и сде-лано на следующем рисунке 11. Границу нового шестиугольного пространства $\bigcirc$ на равных правах представляют вырожденные треугольники первого и второго типов. Вершинам шестиугольника соответствуют вырожденные равнобедренные треугольники с основанием нулевой длины и двумя прилегающими нему углами величиной 0° и 180°.

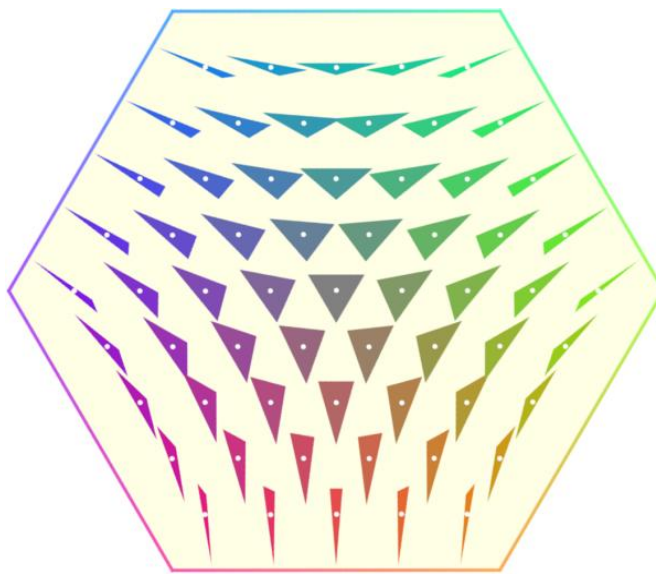

Рис. 11. Пространство треугольников $\bigcirc$

Добавим, что медианное отображение $\mathcal{M}_0: \bigcirc \to \Delta_0$ также допускает непрерывное продолжение на границу $\partial\bigcirc$.

**Литература**

**Балк М.Б, Болтянский В.Г.** (1987). Геометрия масс.

**Behrend K.** (2014). Introduction to algebraic stacks, in *Moduli Spaces*, Ed. by L. Brambila-Paz, P. Newstead, R. P.Thomas, and O. García-Prada, London Mathematical Society Lecture Notes **411**, Cambridge Univ. Press, Cambridge, UK, 1–131.

**Kimberling C.** (2020). Encyclopedia of triangle centers. Available at:
http://faculty.evansville.edu/ck6/encyclopedia/

**Панов П.А.** (2018)**.** О геометрической медиане выпуклых, а также треугольных и других многоугольных областей, Известия Иркутского государственного университета. Серия Математика, 26 , 62–75

**Панов П.А., Савватеев А.В.** (2020). О геометрической медиане и других медианоподобных точках. Экономика и математические методы, 56(3), с. 146–153.

**Панов П.А.** (2020). Визуализация медианного отображения $\mathcal{M}_1$
https://drive.google.com/file/d/1mdSaF_nLsUAKy8xOTWO3kM_emeUa_cxR/view